\documentclass{amsart}
\usepackage{amsfonts,amssymb,amsmath,amsthm}
\usepackage{url}
\usepackage{enumerate}
\usepackage{fullpage}
\usepackage{bm}
\usepackage[usenames,dvipsnames]{color}

\newtheorem{thm}{Theorem}[section]

\theoremstyle{definition}

\numberwithin{equation}{section}

\newcommand{\C}{\mathbb{C}}
\newcommand{\F}{\mathbb{F}}

\newcommand{\Z}{\mathbb{Z}}
\newcommand{\Q}{\mathbb{Q}}
\newcommand{\PP}{\mathbb{P}}
\newcommand{\N}{\mathbb{N}}
\newcommand{\TT}{\mathbb{T}}
\newcommand{\m}{\mathrm{m}}
\newcommand{\re}{\mathop{\mathrm{Re}}}
\newcommand{\im}{\mathop{\mathrm{Im}}} 

\author{Marie-Jos\'e Bertin, Amy Feaver, Jenny Fuselier, Matilde Lal\'{i}n, Michelle Manes}

 \address{Marie-Jos\'e Bertin: Universit\'e Pierre et Marie Curie (Paris 6), Institut de Math\'ematiques, 175 rue
du Chevaleret, 75013 Paris, France}\email{bertin@math.jussieu.fr}

\address{Amy Feaver: Department of Mathematics, University of Colorado at Boulder, Campus Box 395, Boulder, CO 80309, USA}\email{amy.feaver@colorado.edu} 

\address{Jenny Fuselier: Department of Mathematics \& Computer Science, Drawer 31, High Point University, 833 Montlieu Ave.,
High Point, NC 27262, USA}\email{jfuselie@highpoint.edu}  

\address{Matilde Lal\'in:  D\'epartement de math\'ematiques et de statistique,
                                    Universit\'e de Montr\'eal.
                                    CP 6128, succ. Centre-ville.
                                     Montreal, QC H3C 3J7, Canada}\email{mlalin@dms.umontreal.ca}

\address{Michelle Manes: Department of Mathematics, University of Hawaii, 2565 McCarthy Mall, Honolulu, HI 96822, USA}\email{mmanes@math.hawaii.edu}
 
\thanks{This work of ML was partially supported by NSERC Discovery Grant 355412-2008 and FQRNT Subvention \'etablissement
de nouveaux chercheurs 144987.  The work of MM was partially supported by NSF-DMS 1102858.}



\begin{document}
\title{Mahler measure of some singular $K3$-surfaces}

\begin{abstract} We study the Mahler measure of the three-variable Laurent polynomial $x + 1/x + y + 1/y + z + 1/z - k$ 
where $k$ is a parameter. The zeros of this polynomial define (after desingularization) a family of $K3$-surfaces. 
In favorable cases, the $K3$-surface has Picard number~$20$,  and the Mahler measure is related to its $L$-function. 
This was first studied by Marie-Jos\'e Bertin. In this work, we prove several new formulas, extending the earlier work of Bertin. 

\end{abstract}
\maketitle

\section{Introduction}
Given a nonzero Laurent polynomial $P \in \C[x_1^{\pm1}, \dots, x_n^{\pm1}]$, the (logarithmic) Mahler
measure is defined by 
\begin{eqnarray*}
\m(P)\!\! &=& \! \!\int_0^1 \dots \int_0^1 \log \left| P(e^{2 \pi i \theta_1}, \dots ,
e^{2 \pi i \theta_n}) \right|  d \theta_1 \cdots d \theta_n\\
&=&\frac{1}{(2  \pi i)^n} \int_{\TT^n} \log \left| P(x_1, \dots, x_n) \right| \frac{
d x_1}{x_1} \cdots \frac{ d x_n}{x_n},
\end{eqnarray*}
where $\TT^n=\{(x_1,\dots,x_n)\in \C^n\,:\,|x_1|=\dots=|x_n|=1\}$ is the unit $n$-torus.

Jensen's formula relates the Mahler measure of a one-variable polynomial to a very simple formula depending on the roots of the polynomial:
\[\m(P)= \log |a| + \sum_{|r_j|>1} \log |r_j|  \quad \mathrm{for} \quad  P(x) = a \prod_j(x-r_j).\]
This formula shows, in particular, that the Mahler measure of a polynomial with integral coefficients is the logarithm of 
an algebraic number.

The situation for several variable polynomials is very different. There are several formulas for specific polynomials yielding special values of $L$-functions. The first examples were computed by 
Smyth in the 1970s~\cite{S, B1} and give special values of the Riemann zeta function and Dirichlet $L$-series:
 \begin{align*}
 \m(x+y+1) 
 &= \frac{3 \sqrt{3}}{4 \pi} L(\chi_{-3},2) = L'(\chi_{-3}, -1),\\
\m(x+y+z+1) &= \frac{7}{2\pi^2}\zeta(3).
\end{align*}
Then, in the mid 1990s,  Boyd \cite{B2} (after a suggestion of  Deninger) looked at more complicated polynomials such as the family 
\begin{equation}\label{Boyd-eq}
P_k(x,y)=x+\frac{1}{x} + y + \frac{1}{y} - k,
\end{equation} where $k$ is
an integral parameter. 
For most $k$, the zero set $P_k(x,y)=0$ is a genus-one curve which we will denote by $E_{(k)}$. Boyd found several 
numerical formulas of the same shape:
\[\m\left(x+\frac{1}{x} + y + \frac{1}{y} - k\right) \stackrel{?}{=}
s_kL'(E_{(k)},0)  \quad k \in \Z, \quad |k| \not = 0,4,\]
where $s_k$ is a rational number and the question mark means that both sides of the equation are equal to at least 25 decimal places.  In fact, it suffices to consider $k$ natural since one can easily see that $\m(P_k)=\m(P_{-k})$.

In particular, for $k=1$,
\begin{equation} \label{mat-theorem}
\m\left(x +\frac{1}{x}+y+\frac{1}{y}-1\right)=\frac{15}{4 \pi^2}L(E_{15},2)=L'(E_{15},0),
\end{equation}
which was recently proven by Rogers and Zudilin \cite{RZ1}. 

The connection with the $L$-function of the elliptic curve defined by the zeros of the polynomial was explained by Deninger \cite{De} in a very general context and Rodriguez-Villegas \cite{RV} for some specific formulas in terms of Beilinson's conjectures. 
Beilinson's conjectures predict that special values of $L$-functions (coming from an arithmetic-geometric structure) 
are given by certain values of the regulator associated to the structure up to a rational number. In favorable cases, Mahler measure can be related to the regulator. In particular, this allowed Rodriguez-Villegas to prove the formulas for the case where $E$ has complex multiplication, since in this case Beilinson's conjectures are known to be true.

More generally, let $P(x,y)$ be a polynomial in two variables with integer coefficients and suppose that $P$ does not vanish
 on the $2$-torus $\mathbb T^2$.
If $P$ defines an elliptic curve $E$ and the polynomials of the faces $P_F$ of $P$ (defined in terms of the Newton polygon of $P$)
 are cyclotomic (in other words, they have measure zero), then the following relation between $\m(P)$ and the $L$-series of the elliptic curve $E$ is conjectured to hold:
\[
\m(P) \stackrel{?}{=} \frac {qN}{4\pi ^2}L(E,2)=qL'(E,0), 
\]
where $N$ is the conductor of $E$ and $q$ is a rational number.

A natural extension to this connection involves polynomials whose zeros define Calabi-Yau varieties. One-dimensional Calabi-Yau varieties are elliptic curves,  
while $2$-dimensional Calabi-Yau varieties are elliptic $K3$-surfaces. For example, it is natural to consider the family of polynomials resulting from 
adding an extra variable to the equation in  (\ref{Boyd-eq}). Bertin has been pursuing this program \cite{Be06, Be08, Be10} with the families
\[P_k(x,y,z)=x+\frac{1}{x} + y + \frac{1}{y} + z +\frac{1}{z} - k.\]  
\begin{equation*}
Q_k(x,y,z)=x+\frac {1}{x}+y+\frac {1}{y}+z+\frac {1}{z} +xy+\frac {1}{xy}+zy+\frac {1}{zy}+xyz+\frac {1}{xyz}-k.
\end{equation*}

Relating these examples back to the elliptic curve case, one may ask for a natural condition on the faces of the Newton polytope for the polynomials $P_k$ in order to expect relationships between $\m(P_k)$ and the $L$-series of the associated surface.  The situation is more complicated than in the elliptic curve case, since the  faces in the above examples have nonzero Mahler measure.  This question remains open.


The first step in Bertin's work is to generalize Rodriguez-Villegas's expression of the Mahler measure in terms of 
Eisenstein--Kronecker series for these two families of polynomials defining $K3$-surfaces.
For example, in \cite{Be06} Bertin proves 
\[\begin{aligned}  \m(P_k)=&\frac {\im \tau}{8 \pi ^3}\sum_{j\in \{1,2,3,6\}} \sum'_{m,n}(-1)^j4j^2\left(2\re
\frac {1}{(jm\tau+n)^3(jm\bar{\tau}+n)}+\frac
{1}{(jm\tau+n)^2(jm\bar{\tau}+n)^2}\right).  \end{aligned} \]
Here
$k=w+\frac {1}{w}$ and 
\[
w=\left(\frac {\eta (\tau) \eta (6\tau)}{\eta
(2\tau)\eta(3\tau)}\right)^6=q^{1/2}-6q^{3/2}+15q^{5/2}-20q^{7/2}+\cdots
\]
where $\eta$ denotes the Dedekind eta function. 

For exceptional values of $k$, the corresponding $K3$-surface $Y_k$ is singular (or extremal) and $\tau$ is imaginary quadratic. The  Eisenstein--Kronecker series
can be split into two sums, one with the $\re
\frac {1}{(jm\tau+n)^3(jm\bar{\tau}+n)}$ terms and the other with the $\frac
{1}{(jm\tau+n)^2(jm\bar{\tau}+n)^2}$ terms. The first one is related to the $L$-series of the surface, while the second one is either zero or may be expressed in terms of a Dirichlet series related to the Mahler measure of the 
$2$-dimensional faces of the Newton polytope of the polynomial $P_k$.


Bertin obtained
\begin{eqnarray*}
\m(P_0)&=&d_3:=\frac{3\sqrt{3}}{4\pi}L(\chi_{-3},2),\\
\m(P_2)&=&4\frac{|\det \TT(Y_2)|^{3/2}}{4\pi^3}L(\TT(Y_2),3)=4\cdot\frac {8\sqrt{8}}{4\pi^3}L(g_8,3), \text{ and}\\
\m(P_{10})&=& \frac{4}{9}\frac{|\det \TT(Y_{10})|^{3/2}}{4\pi^3}L(\TT(Y_{10}),3)+2d_3 = \frac{4}{9}\cdot\frac {72 \sqrt{72}}{4\pi^3} L(g_8,3)+2d_3,\\
\end{eqnarray*}
where $Y_k$ denotes the $K3$-surface associated to the zero set $P_k(x,y,z)=0$, $\TT$ denotes its transcendental lattice, and $L(g_N,3)$ the $L$-series at $s=3$ of a modular form of weight $3$ and level $N$.


In this note, we continue the work of Bertin and prove
\begin{eqnarray*}
 \m(P_3)&=&2\frac{|\det \TT(Y_{3})|^{3/2}}{4\pi^3}L(\TT(Y_3),3)=2\cdot \frac{15\sqrt{15}}{4\pi^3}L(g_{15},3),\\
\m(P_6)&=&2\frac{|\det \TT(Y_{6})|^{3/2}}{4\pi^3}L(\TT(Y_6),3)=2\cdot \frac{24\sqrt{24}}{4\pi^3}L(g_{24},3), \text{ and}\\
\m(P_{18})&=&\frac{1}{5}\frac{|\det \TT(Y_{18})|^{3/2}}{4\pi^3}L(\TT(Y_{18}),3)+\frac{14}{5}d_3 = \frac{1}{5}\cdot\frac{120\sqrt{120}}{4\pi^3}L(g_{120},3)+\frac{14}{5}d_3.
\end{eqnarray*}

The case with $k=18$ is particularly difficult because the corresponding $K3$-surface has 
an infinite section that is defined over a quadratic field rather than being defined over $\Q$. 
The method we use to find this infinite section should be useful in other cases.

\section{Background on $K3$-surfaces}

A $K3$-surface is a complete smooth surface $Y$ that is simply connected and admits a unique (up to scalars) holomorphic 2-form $\omega$.
 We list here some useful facts  about $K3$-surfaces along with notation that will be used throughout. See \cite{Yu04} for  general results about Calabi-Yau manifolds including $K3$-surfaces.

\begin{itemize}

\item  $H_2(Y, \Z)$ is a free group of rank 22.

\item The Picard group $\mathrm{Pic}(Y)\subset H_2(Y, \Z)$ is the group of divisors modulo linear equivalence, parametrized by algebraic cycles:
\[\mathrm{Pic}(Y) \cong \Z^{\rho(Y)}.\]
The exponent $\rho(Y)$ is called the Picard number, and over a field of characteristic~$0$ it satisfies
\[1 \leq \rho(Y)\leq 20.\]
If $\rho(Y)=20$, we say that the $K3$-surface is \emph{singular.}

\item The transcendental lattice is defined by  
\[\TT(Y)=(\mathrm{Pic}(Y))^\bot.\]

 \item Let $\{\gamma_1, \dots, \gamma_{22}\}$ be a $\Z$-basis for $H_2(Y, \Z)$. Then

\[
\int_\gamma \omega 
=
\left\{\begin{array}{ll} 0 & \gamma\in\mathrm{Pic}(Y), \\
\text{period of }  Y& \gamma\in \TT(Y).
\end{array}\right.\]

\end{itemize}

\subsection{$L$-functions}

Let $Y$ be a surface. 
The zeta function is defined by 
\[Z(Y,u)=\exp \left(\sum_{n=1}^\infty N_n(Y) \frac{u^n}{n} \right), \qquad |u|<\frac{1}{p},\]
where $N_n(Y)$ denotes the number of points on $Y$ in $\F_{p^n}.$

If $Y$ is a $K3$-surface defined over $\Q$, then $Y$ gives a $K3$-surface over $\F_p$ for almost all $p$ and
\[Z(Y,u)=\frac{1}{(1-u)(1-p^2u)P_2(u)},\]
where $\deg P_2(u)=22$.  In fact,
\[P_2(u)=Q_p(u)R_p(u),\]
where the polynomial $R_p(u)$ comes from the algebraic cycles and $Q_p(u)$ comes from the transcendental cycles. 
Hence, for a singular 
$K3$-surface, $\deg Q_p=2$ and $\deg R_p=20$. 

Finally, we will work with the part of the $L$-function of $Y$ coming from the transcendental lattice, which is given by
\[
L(\TT(Y),s)=(*)\prod_{p\,\mathrm{good}}\frac{1}{Q_p(p^{-s})}=\sum_{n=1}^{\infty} \frac{A_n}{n^s},
\]
where  $(*)$ represents finite factors coming from the primes of bad reduction. 

\subsection{Elliptic surfaces}
An elliptic surface $Y$ over $\PP^1$ is a smooth projective surface $Y$ with an elliptic fibration, i.e., a surjective morphism
\[\Phi: Y \rightarrow \PP^1\]
such that almost all of the fibers are smooth curves of genus 1 and no fiber contains an exceptional curve of the first kind (with self-intersection $-1$). Here we list some facts
about elliptic surfaces. See \cite{SS10} for a comprehensive reference containing these results. 

The group of global sections of the elliptic surface is called the Mordell-Weil group and can be naturally identified with the group of points of the generic fiber. 
Its rank $r$ can be found from the formula
\begin{equation} \label{shioda1}
\rho(Y)=r +2 +\sum_{\nu=1}^h(m_\nu -1)
\end{equation}
due to Shioda \cite{Sh90}. Here $m_\nu$ denotes number of irreducible components of the corresponding singular fiber and $h$ is the number of singular fibers.

Global sections can be also thought as part of the  N\'eron-Severi group $\mathrm{NS}(Y)$ given by the divisors modulo algebraic equivalence. 
It is finitely generated and torsion-free. Intersection of divisors yields a bilinear pairing 
which gives  $\mathrm{NS}(Y)$ the structure of an integral lattice. 

The trivial lattice $\mathrm{T}(Y)$ is the subgroup of  $\mathrm{NS}(Y)$ generated by the zero section and the fiber components. Its determinant is given by 
\begin{equation}\label{mutantshioda2}
\det \mathrm{T}(Y) = \prod_{\nu=1}^hm_\nu^{(1)},
\end{equation}
where $m_\nu^{(1)}$ indicates the number of single components of the corresponding singular fiber.  (See~\cite[p. 17]{Sh90}.)
One has that the Mordell-Weil group is isomorphic to 
$\mathrm{NS}(Y)/\mathrm{T}(Y)$.

The Mordell-Weil group can also be given a lattice structure $\mathrm{MWL}(Y)$. Then
\begin{equation}\label{totallycrazydeterminant}
\det \mathrm{NS}(Y)=(-1)^{r} \frac{\det \mathrm{T}(Y) \det \mathrm{MWL}(Y)}{|E_\mathrm{tors}|^2},
\end{equation}
where $E$ is the generic fiber.
The bilinear pairing induced by intersection can be used to construct a height that satisfies
\begin{equation}\label{heightdef}
h(P)=2 \chi(Y) +2 (\overline{P}\cdot \overline{O})-\sum_{\nu}\mathrm{contr}_\nu(P),
\end{equation}
where $ \chi(Y)$ is the arithmetic genus ($\chi(Y)=2$ for $K3$-surfaces), $\overline{P}\cdot \overline{O}\geq 0$, and the (always nonnegative) correction terms $\mathrm{contr}_\nu(P)$ measure how $P$ intersects the components of the singular fiber over $\nu$.  This height is the canonical height that one obtains
by thinking about the elliptic surface as an elliptic curve over a function field~\cite{Sh90}.



\subsection{A particular family of $K3$-surfaces} 
In this note, we  consider the family of polynomials 
\[
P_k(x,y,z)=x+\frac1x+y+\frac1y+z+\frac1z-k. 
\]
The  desingularization of $P_k = 0$ results in a $K3$-hypersurface $Y_{k}$.
We homogenize the numerator of $P_k$:
\[
x^2yz+xy^2z+xyz^2+t^2(xy+xz+yz)-kxyzt,
\]
and then get an elliptic fibration by setting $t=s(x+y+z)$. 
\begin{equation}
Y_k:
s^2(x+y)(x+z)(y+z)+(s^2-ks+1)xyz=0.
\label{eq:ell fib}
\end{equation}
To study the components of the singular fibers, one expresses the $K3$-surface $Y_k$ as a double covering of a well-known rational elliptic surface 
given by Beauville \cite{Bea82}
\begin{equation}\label{beauville}
(x+y)(x+z)(y+z)+uxyz=0.
\end{equation}

By analyzing the structure of the singular fibers, we can compute the rank of the group of sections $r$.  
In the case of Beauville's surface, the singular fibers  are given by
\[\begin{array}{ll}
u=\infty & I_6,\\
u=0 & I_3,\\
u=1 & I_2, \text{ and}\\
u=-8 & I_1.   
  \end{array}\]

To conclude this section, we summarize some results from Peters and Stienstra \cite{PS89} on this family of $K3$-surfaces.  For generic $k$, the Picard number is $\rho(Y_k)=19$.  We focus on the singular $K3$-surfaces --- that is, on $k$ values for which   $\rho(Y_k)=20$.  The transcendental lattice $\TT$ of the general family $Y_k$ has a Gram matrix 
of the form
\begin{equation}
\left(\begin{array}{ccc} 0 & 0 & 1\\ 0 & 12 & 0 \\ 1 & 0 &0        \end{array}
 \right).\label{eqn:gram mat}
 \end{equation}
Having Picard number $\rho=20$ is equivalent to having a relation between the generic basis $\{\gamma_1,\gamma_2,\gamma_3\}$ of transcendental periods; that is, 
\begin{equation}
p\gamma_1+q\gamma_2+r\gamma_3
\label{eqn:becomealg}
\end{equation}
becomes algebraic for some choice of $p,q,r$.

 Now, let $k=w+\frac{1}{w}$. Then $w$ can be represented as a modular function:
\[
w=\left(\frac{\eta(\tau)\eta(6\tau)}{\eta(2\tau)\eta(3\tau)}\right)^6, \quad \eta(\tau)=e^\frac{\pi i \tau}{12} \prod_{n\geq 1}(1-e^{2\pi i n \tau}),\quad \tau \in \mathbb{H},
\]
where $\mathbb{H}$ denotes the upper half-plane.  Furthermore, a period is algebraic precisely when it is orthogonal to 
$\gamma_1 + \tau \gamma_2 - 6\tau^2 \gamma_3$.
Combining these facts yields a quadratic equation for $\tau$: 
\begin{equation}
-6p\tau^2+12q\tau+r=0.
\label{eqn:taueq}
\end{equation} 
Thus to find $k$-values such that $Y_k$ is a singular $K3$-surface,  
 we look for $k$ values yielding an imaginary quadratic $\tau$.   Here are a few such values:
\begin{center}
\begin{tabular}{ |c | c|c|c|c|c|c|}\hline\label{table: k tau}
$k$ & $0$ & $2$ & $3$ & $6$ & $10$ & $18$\\
\hline
$\tau$ & $\frac{-3+\sqrt{-3}}{6}$ &$\frac{-2+\sqrt{-2}}{6}$& $\frac{-3+\sqrt{-15}}{12}$ & $\frac{1}{\sqrt{-6}}$ & $\frac{1}{\sqrt{-2}}$&$\sqrt{\frac{-5}{6}}$ \\ 
\hline
\end{tabular}
\end{center}
Given $\tau$, one may find the parameters $p$, $q$, and $r$, and then find the discriminant of $\TT$ up to squares by taking the determinant of the resulting Gram matrix.  See Section~\ref{sec:lat and rank} for details in the cases where $k=3$, $k=6$, and $k=18$.

\section{Main results and the general strategy for the proof}

\begin{thm} We have the following formulas:
\begin{eqnarray*}
\m(P_3)&=&\frac{15\sqrt{15}}{2\pi^3}L(g_{15},3)=2\frac{|\det \TT(Y_{3})|^{3/2}}{4\pi^3}L(\TT(Y_3),3),\\
\m(P_6)&=&\frac{24\sqrt{24}}{2\pi^3}L(g_{24},3)=2\frac{|\det \TT(Y_{6})|^{3/2}}{4\pi^3}L(\TT(Y_6),3), \text{ and}\\
\small
\m(P_{18})&=&\frac{120\sqrt{120}}{20\pi^3}L(g_{120},3)+\frac{14}{5}d_3=\frac{1}{5}\frac{|\det \TT(Y_{18})|^{3/2}}{4\pi^3}L(\TT(Y_{18}),3)+\frac{14}{5}d_3,
\end{eqnarray*}
where $Y_{k}$ is the $K3$-hypersurface defined by the zeros of $P_k(x,y,z)$, $\TT(Y_k)$ is its transcendental lattice, and $g_N$
is a CM modular form of level $N$. 
\end{thm}

The strategy for proving these formulas is as follows:

\begin{itemize}

\item Understand the transcendental lattice and the group of sections. 

\item Relate the Mahler measure $\m(P_k)$ to the $L$-function of a modular form.

\item Relate the $L$-function of the surface $Y_k$ to the $L$-function of that same modular form.

\end{itemize}

\section{The Transcendental Lattice and the Rank}\label{sec:lat and rank}
We will prove the following:

\begin{itemize}

 \item For $k=6$, $|\det \TT({Y_6})|=24$, $\text{rank}=0$.

\item For $k=3$, $|\det \TT({Y_3})|=15$, $\text{rank}=1$.

\item For $k=18$, $|\det \TT({Y_{18}})|=120$, $\text{rank}=1$.
\end{itemize}

\subsection{The transcendental lattice and the rank for $Y_6$}\label{sec: lat and rank 6}
When $k=6$, we see from the table on page~\pageref{table: k tau} that $\tau=\frac{1}{\sqrt{-6}}$. Thus, it satisfies the equation $-6\tau^2-1=0$, so in equation~\eqref{eqn:taueq} we take $p=1$, $q=0$, and $r=-1$.  By equation~\eqref{eqn:becomealg}, the   
vector $\gamma_1-\gamma_3$ becomes algebraic over $Y_6$.   That is,
$v=\gamma_1 - \gamma_3 \in \mathrm{Pic}(Y_6)$.  

To find the transcendental lattice, we use the Gram matrix~\eqref{eqn:gram mat} to find  vectors  orthogonal to~$v$. A simple computation yields: $\{\gamma_2, \gamma_1+\gamma_3\}$; hence these span a sublattice of $\TT$.
We again use~\eqref{eqn:gram mat}, this time to find the  Gram matrix for the space spanned by these two vectors:
\[\left( \begin{array}{cc}12 & 0 \\ 0 & 2\end{array}\right).\]
Thus the discriminant of $\TT$, up to a square, is equal to 24.   It remains to decide if it is~$6$ or~$24$.

Equation~\eqref{eq:ell fib} expresses $Y_6$ as a double-covering of the Beauville surface~\eqref{beauville},  with $u = (s^2-6s+1)/s^2$.
\[Y_6:
s^2(x+y)(x+z)(y+z)+(s^2-6s+1)xyz=0.\]
 Since we know the singular fibers of the Beauville surface, we easily find the singular fibers of~$Y_6$: 
\[\begin{array}{lcl}
s=0 & I_{12} & \mbox{double over } u=\infty,\\
s=\alpha & I_3 & \mbox{over } u=0,\\
s=\beta & I_3 & \mbox{over } u=0,\\
s=\frac{1}{6} & I_2 & \mbox{over } u=1,\\
s=\infty & I_2 & \mbox{over } u=1, \text{ and}\\  
s=\frac{1}{3} & I_2 & \mbox{double over } u=-8.
  \end{array}\]
(Here $\alpha$ and $ \beta$ are the two distinct roots of $s^2-6s+1=0$.)

Applying  Shioda's formula~\eqref{shioda1}, we have
\[20=r+2+(12-1)+(3-1)+(3-1)+(2-1)+(2-1)+(2-1)=r+20,\]
so the rank of the group of sections is~$0$. 
A Weierstrass form is given by 
\[
y^2+(s^2-6s+1)xy=x(x-s^4)(x+s^2-6s^3).
\]
 We can compute the torsion group directly. A point of order 6 is given by 
\[
\left(s^2(6s-1),0\right)
\]
and  the only point of order 2 is $(0,0)$.

Applying formula~\eqref{totallycrazydeterminant}, we have
\[
\left|\det \TT(Y_6)\right|
=
|\det \mathrm{NS}(Y_6)|
=
\frac{12\cdot 3 \cdot 3 \cdot 2 \cdot 2 \cdot 2}{|E_\mathrm{tors}|^2}
=\frac{2^5\cdot 3^3}{|E_\mathrm{tors}|^2}.\]
This means that  either  ${|E_\mathrm{tors}|}=6$ and $|\det \TT_{Y_6}|=24$, or ${|E_\mathrm{tors}|}=12$  and $|\det \TT_{Y_6}|=6$.
By the work of Miranda and Persson \cite{MP}, ${|E_\mathrm{tors}|}=12$  implies that the torsion is given by $\Z/6\Z \times \Z/2\Z$ which is not possible since 
there is only one point of order 2. Therefore,  ${|E_\mathrm{tors}|}=6$ and
\[|\det \TT(Y_6)|=24.\]

\subsection{The transcendental lattice and the rank for $Y_3$}
In this case we have  $\tau=\frac{-3+\sqrt{-15}}{12}$ (see the table on page~\pageref{table: k tau}), which satisfies the quadratic equation $-6\cdot 4 \tau^2-12\tau-4=0$.  So in equation~\eqref{eqn:taueq} we take
$p=4$, $q=-1$, and $r=-4$.  By equation~\eqref{eqn:becomealg}, $v=4\gamma_1-\gamma_2-4\gamma_3 \in \mathrm{Pic}(Y_3)$. Using the Gram matrix~\eqref{eqn:gram mat}, we find that  $\{\gamma_1+\gamma_3, \gamma_2+3\gamma_3\}$ generate a sublattice of $\TT$, and their Gram matrix is:
\[\left( \begin{array}{cc}2 & 3 \\ 3 & 12\end{array}\right).\]
Since  the determinant of this matrix is square-free,  we conclude that $|\det \TT(Y_3) | = 15$. 

The equation
\[s^2(x+y)(x+z)(y+z)+(s^2-3s+1)xyz=0,\]
expresses $Y_3$ as a double-covering of the Beauville surface~\eqref{beauville} with $u =  (s^2-3s+1)/s^2$.
In this case, the singular fibers are:
\[\begin{array}{lcl}
s=0 & I_{12} & \mbox{double over } u=\infty,\\
s=\alpha_1 & I_3 & \mbox{over } u=0,\\
s=\beta_1 & I_3 & \mbox{over } u=0,\\
s=\frac{1}{3} & I_2 & \mbox{over } u=1,\\
s=\infty & I_2 & \mbox{over } u=1,\\  
s=\alpha_2 & I_1 & \mbox{over } u=-8, \text{ and}\\
s=\beta_2 & I_1 & \mbox{over } u=-8.\\
  \end{array}\]
Here, $\alpha_1, \beta_1$ are the two distinct roots of $s^2-3s+1=0$, and $\alpha_2, \beta_2$ are the  roots
of $9s^2-3s+1=0$.

By Shioda's formula \eqref{shioda1}, the rank is~$1$.  A Weierstrass model around infinity is given by:
\[
y^2+(\sigma^2-3\sigma+1)xy=x(x-1)(x+\sigma^2-3\sigma)=x^3+(\sigma^2-3\sigma-1)x^2+(-\sigma^2+3\sigma)x.
\]

With the aid of Pari/gp or Sage~\cite{PARI, St11} we find a point $\rho_6$ of order 6. Indeed,
\begin{eqnarray*}
\rho_6&=&\left(-\sigma(\sigma-3),\sigma(\sigma-3)(\sigma^2-3\sigma+1)\right),\\
2\rho_6&=&\left(1,-\sigma^2+3\sigma-1\right),\\
3\rho_6&=&\left(0,0\right),\\
4\rho_6&=&\left(1,0\right), \text{ and}\\
5\rho_6&=&\left(-\sigma^2+3\sigma,0\right).
\end{eqnarray*}
By the work of Miranda and Persson~\cite{MP}, since the rank is~$1$ and $\chi=2$,  the torsion must have order~$6$, and therefore
it must be generated by $\rho_6$. 

With the aid of Pari/gp or Sage we also find the following point in each fiber:
\[
\left(-(\sigma-3)(\sigma-1)^2, (\sigma-3) (\sigma-2) (\sigma-1) (\sigma^2-3 \sigma+1)\right).
\]
Since this point is not generically among the torsion points of each fiber, it must give an infinite section, 
which is in particular defined over $\Q$. In fact, this point is a generator of the infinite section, but we do not need this fact for our computation.

\subsection{The transcendental lattice and the rank for $Y_{18}$} \label{infi18}
When $k=18$, the table shows $\tau=\sqrt{\frac{-5}{6}}$, which satisfies $-6\tau^2-5=0$.  Take $p=1$, $q=0$, and  $r=-5$ in equation~\eqref{eqn:taueq}, so  $v=\gamma_1-5\gamma_3 \in \mathrm{Pic}(Y_{18})$.  The vectors  $\{\gamma_2, \gamma_1+5\gamma_3\}$ are orthogonal to $v$, and the corresponding Gram matrix is
\begin{equation}\label{eqn: gram 18}
\left( \begin{array}{cc}12 & 0 \\ 0 & 10\end{array}\right).
\end{equation}
The determinant of this matrix is 120, so the discriminant of the transcendental lattice is either 30 or 120. 

The double-cover of the Beauville surface is given by:
\[
Y_{18}: s^2(x+y)(x+z)(y+z)+(s^2-18s+1)xyz=0,\]
where $u = (s^2-18s+1)/s^2$.
The singular fibers are
\[\begin{array}{lcl}
s=0 & I_{12} & \mbox{double over } u=\infty,\\
s=\alpha_1 & I_3 & \mbox{over } u=0,\\
s=\beta_1 & I_3 & \mbox{over } u=0,\\
s=\frac{1}{18} & I_2 & \mbox{over } u=1,\\
s=\infty & I_2 & \mbox{over } u=1,\\  
s=\alpha_2 & I_1 & \mbox{over } u=-8, \text{ and}\\
s=\beta_2 & I_1 & \mbox{over } u=-8.\\
  \end{array}\]
Here $\alpha_1, \beta_1$ are the two distinct roots of $s^2-18s+1=0$, and
 $\alpha_2, \beta_2$ are the  roots of $9s^2-18s+1=0$.

From Shioda's formula~\eqref{shioda1}, we see that the rank is~$1$.
A Weierstrass model around infinity is given by 
\begin{equation}
y^2+(\sigma^2-18\sigma+1)xy=x(x-1)(x+\sigma^2-18\sigma)=x^3 + (\sigma^2 - 18\sigma - 1)x^2 + (-\sigma^2 + 18\sigma)x.
\label{eqn:Y18 model}
\end{equation}
With the aid of  Pari/gp or Sage~\cite{PARI, St11}, we find a point $\rho_6$ of order 6. Indeed,
\begin{eqnarray*}
\rho_6&=&\left(-\sigma(\sigma-18),  \sigma (\sigma-18) \left(\sigma^2-18 \sigma+1\right)\right),\\
2\rho_6&=&\left(1,-\sigma^2+18\sigma-1\right),\\
3\rho_6&=&\left(0,0\right),\\
4\rho_6&=&\left(1,0\right), \text{ and}\\
5\rho_6&=&\left(-\sigma^2+18\sigma,0\right).\\
\end{eqnarray*}
Again by the work of Miranda and Persson \cite{MP},  $r=1$ and $\chi=2$ implies that the torsion must have order 6, and hence must be generated by $\rho_6$.

If $P$ is a generator of the infinite part of the group of sections, then $\det \mathrm{MWL}(Y_{18})=h(P)$. Applying formulas~\eqref{mutantshioda2} and~\eqref{totallycrazydeterminant}, we have
\begin{equation}\label{heightformula}
|\det\TT(Y_{18})|= |\det \mathrm{NS}(Y_{18})|= \frac{12\cdot 3^2\cdot 2^2 h(P)}{6^2}=12h(P).
\end{equation}
By the remark following ~\eqref{eqn: gram 18}, $|\det\TT(Y_{18})| = 30$ or 120.  Hence either $|\det\TT(Y_{18})|=30$ and $h(P)=5/2$ or $|\det\TT(Y_{18})|=120$ and $h(P)=10$.

Finding the infinite section for $Y_{18}$ is more difficult than for $Y_3$ because the infinite section is not defined over $\Q$.  Details of the method used to find the infinite section, prove that we have a generator, and compute its height are in Section~\ref{sec:Inf Secs}.
 The outcome of the computations is  a generator   $p_\sigma$ defined over $\Q(\sqrt{-3})$ satisfying $h(p_\sigma)=10$; hence
\[|\det \TT(Y_{18}) | = 120.\]


\section{Relating the Mahler Measure  to a newform}
The main ingredient we use to relate Mahler mesure to newforms is the following result from~\cite{Be06}.

\begin{thm} [Bertin]\label{thm: bertin thm sum}
 Let $k= w +\frac{1}{w}$ with
\[w=\left(\frac{\eta(\tau)\eta(6\tau)}{\eta(2\tau)\eta(3\tau)}\right)^6, \quad
 \eta(\tau)=e^\frac{\pi i \tau}{12} \prod_{n\geq 1}(1-e^{2\pi i n \tau}).\]
Then
\begin{eqnarray*}
\hspace{-0.1cm}\m(P_k)\hspace{-0.3cm}&=&\hspace{-0.3cm}\frac{\im \tau}{8 \pi^3}\left[\sum_{m,n}\left(-4\left(2\re \frac{1}{(m\tau+n)^3(m\bar{\tau}+n)}+\frac{1}{(m\tau+n)^2(m\bar{\tau}+n)^2} \right) \right.\right.\\
&&\qquad \qquad
+16 \left(2\re \frac{1}{(2m\tau+n)^3(2m\bar{\tau}+n)}+\frac{1}{(2m\tau+n)^2(2m\bar{\tau}+n)^2} \right)\\
&&\qquad \qquad
-36\left(2\re \frac{1}{(3m\tau+n)^3(3m\bar{\tau}+n)}+\frac{1}{(3m\tau+n)^2(3m\bar{\tau}+n)^2} \right)\\
&&\qquad \qquad
+144\left.\left(2\re \frac{1}{(6m\tau+n)^3(6m\bar{\tau}+n)}+\frac{1}{(6m\tau+n)^2(6m\bar{\tau}+n)^2} \right)\right].
\end{eqnarray*}
\end{thm}

The evaluation of the Eisenstein--Kronecker series often leads  to Hecke $L$-functions.  Let $K$ be an imaginary quadratic number field and $\mathfrak{m}$ be an ideal of $\mathcal{O}_K$. A Hecke character of
$K$ modulo $\mathfrak{m}$ with $\infty$-type $\ell$ is a homomorphism $\phi$ on the group of fractional ideals of $K$ which are prime
to $\mathfrak{m}$ such that for all $\alpha \in K^*$	 with $ \alpha \equiv 1 \mod \mathfrak{m}$,
\[\phi((\alpha))=\alpha^\ell  .\] 
The ideal $\mathfrak{m}$ is called the conductor of $\phi$ if it is minimal in the following sense: if $\phi$ is defined modulo  $\mathfrak{m}'$, then $\mathfrak{m} | \mathfrak{m}'$.

Let \[L(\phi, s)= \sum_{\mathfrak{a} \, \mathrm{integral}}\frac{\phi(\mathfrak{a})}{N(\mathfrak{a})^s} =\sum_{cl(\mathfrak{a})}\frac{\phi(\mathfrak{a})}{N(\mathfrak{a})^{2-s}}\frac{1}{2}\sum_{\lambda \in \mathfrak{a}}' \frac{\bar{\lambda}^2}{(\lambda\bar{\lambda})^s}.\]
The Mellin transform gives a Hecke eigenform:
\[f_\phi=\sum_{n\in \N} a_nq^n =\sum_{\mathfrak{a} \, \mathrm{integral}} \phi(\mathfrak{a}) q^{N(\mathfrak{a})}.\]
A theorem of Hecke and Shimura implies that $f_\phi$ has  weight $\ell +1$ and level $\Delta_KN(\mathfrak{m})$. If $\ell$ is even,
\[f_\phi\in S_{\ell+1}(\Gamma_0(\Delta_KN(\mathfrak{m})), \chi_K)\]
where $-\Delta_K$ is the discriminant of the field, and $\chi_K$ is its quadratic character.

A newform $f=\sum a_n q^n \in S_k(\Gamma_1(N))$ is said to have complex multiplication (CM) by a Dirichlet character~$\phi$ if 
$f=f \otimes \phi$, where
\[f\otimes \phi= \sum_{n\in \N}\phi(n) a_n q^n.\] 
By a result of Ribet, a newform has CM by a quadratic field $K$ iff it comes from a Hecke character of $K$. In particular, $K$ 
is imaginary and unique. 
Sch\"utt \cite{Sch08} proves that there are only finitely many CM newforms with rational coefficients for certain fixed weights (including 3) up to twisting, and he gives a comprehensive table for these.

\subsection{The relation with a newform for $P_6$}
From Theorem~\ref{thm: bertin thm sum},
\[
\m(P_6)=\frac{24\sqrt{6}}{\pi^3}\left(\frac{1}{2} \sum_{m,k}' \left(\frac{m^2-6k^2}{(m^2+6k^2)^3} +\frac{3k^2-2m^2}{(3k^2+2m^2)^3} \right)\right).
\]

This summation can be viewed (see \cite{Be06}) as a Hecke $L$-series on the field $\Q(\sqrt {-6})$. This field has discriminant~$-24$ and class number~$2$, with the nontrivial class represented by $(2,\sqrt{-6})$. That is, we have 
\[\m(P_6)= \frac{24\sqrt{6}}{\pi^3} L_{\Q(\sqrt{-6})}(\phi, 3),
\text{ where } \phi(2,\sqrt{-6})=-2. 
\]

 By the results of Hecke and Shimura, we look for  a correspondence to a (quadratic) twist of  a newform of weight~$3$ and level~$24$.
According to Sch\"utt's table~\cite{Sch08}, there is only one newform (up to twisting) of weight~$3$ and level~$24$.
The twist must be of the form $\left( \frac d p \right ) $ for $d$ dividing $24$, and we can compute the twist exactly by 
comparing the first few coefficients, as shown in the following table.

\bigskip

\begin{center}

\label{table6}
\begin{tabular}{|c || c| c| c| c| c| c| c| c| c| c| c| c|} \hline
$a_p$ & 2 & 3 & 5 & 7 & 11 & 13 & 17 & 19 & 23 & 29 & 31\\\hline \hline
newform of level 24 & 2 & $-3$ & $-2$ & $-10$ & 10 & 0 &0 & 0 & 0 & $-50$ & 38 \\ \hline
coef. of $L_{\Q(\sqrt{-6})}(\phi, s)$ & $-2$ & 3 & 2 & $-10$ & $-10$ & 0 &0 & 0 & 0 & 50 & 38 \\ \hline
\end{tabular}

\end{center}

\bigskip

We find that the twist is given by $ \left( \frac{-3} p \right)$. Therefore, 
\begin{equation}
\m(P_6)
=\frac{24\sqrt{6}}{\pi^3}L\left(f_{24}\otimes \left(\frac{-3}{\cdot}\right),3\right).
\label{eqn: P6 mahler L}
\end{equation}


\subsection{The relation with a newform for  $P_3$}
This case was also considered in \cite{Be06} as a Hecke $L$-series on the field $\Q(\sqrt {-15})$. This field has discriminant~$-15$ and class number~$2$, with the nontrivial class  represented by 
$\left(2,\frac{1+\sqrt{-15}}{2}\right)$.
\begin{eqnarray*}
\m(P_3)&=&\frac{15\sqrt{15}}{2\pi^3}\left(\frac{1}{4}\sum_{m,k}' \left(\frac{2m^2+2mk-7k^2}{\left(m^2+mk+4k^2\right)^3}-\frac{m^2+8mk+k^3}{\left(2m^2+mk+2k^2\right)^3} \right)\right)\\
&=& \frac{15\sqrt{15}}{2\pi^3}L_{\Q(\sqrt{-15})}(\phi,3),
\end{eqnarray*}
where $\phi\left(2,\frac{1+\sqrt{-15}}{2}\right)=-2$. 

There is only one newform of level~$15$ and weight~$3$ in Sch\"utt's table.  We compare the first few coefficients.

\bigskip

\begin{center}
\begin{tabular}{|c || c| c| c| c| c| c| c| c| c| c| c| c|}\hline
$a_p$ & 2 & 3 & 5 & 7 & 11 & 13 & 17 & 19 & 23 & 29 & 31\\\hline \hline
newform of level 15 & $-1$ & 3 & $-5$ & 0 & 0 & 0 & 14 & $-22$ & $-34$ & 0 & 2 \\ \hline
coef. of $L_{\Q(\sqrt{-15})}(\phi, s)$ & $-1$ & 3 & $-5$ & 0 & 0 & 0 & 14 & $-22$ & $-34$ & 0 & 2 \\ \hline
\end{tabular}
\end{center}
\bigskip

Therefore,
\begin{equation}
\m(P_3)= \frac{15\sqrt{15}}{2\pi^3}L\left(f_{15},3\right).
\label{eqn: P3 mahler L}
\end{equation}

\subsection{The relation with a newform for  $P_{18}$}

After some algebraic manipulation, one can find a Hecke series in $\Q(\sqrt{-30})$. This field has discriminant~$-120$  and class number 4, with the class group generated by $(2,\sqrt{-30})$ and $(3,\sqrt{-30})$. We have
\begin{align*}
\m(P_{18})
& = \frac{6 \sqrt{120}} { \pi^3}\left( \frac{1}{2}\sum_{m,k}' \left (\frac{5m^2-6k^2}{(5m^2+6k^2)^3}-\frac{10m^2-3k^2}{(10m^2+3k^2)^3}+\frac{15m^2-2k^2}{(15m^2+2k^2)^3}-\frac{30m^2-k^2}{(30m^2+k^2)^3}\right)\right)\\
&\quad +\frac{3\sqrt{30}}{\pi^3}\sum_{m,k}'\left (-\frac{1}{(5m^2+6k^2)^2}+\frac{1}{(10m^2+3k^2)^2}-\frac{1}{(15m^2+2k^2)^2}+\frac{1}{(30m^2+k^2)^2} \right)\\
& = \frac{6\sqrt{120}}{\pi^3}L_{\Q(\sqrt{-30})}(\phi,3)+\frac{14}{5}d_3,
\end{align*}
where $\phi(2,\sqrt{-30})=-2$ and $\phi(3,\sqrt{-30})=3$.
The equality for the term $\frac{14}{5}d_3$ was proved by Bertin \cite{Be11} by examining identities of certain
 Epstein zeta functions. 
 
There is only one newform of weight~$3$ and level~$120$ in Sch\"utt's table. 

\begin{center}
\bigskip

\begin{tabular}{|c || c| c| c| c| c| c| c| c| c| c| c| c|}\hline
$a_p$ & 2 & 3 & 5 & 7 & 11 & 13 & 17 & 19 & 23 & 29 & 31\\\hline \hline
newform of level 120 & 2 & 3 & $-5$ & 0 & 2 & $-14$ & $-26$ & 0 & $-14$ & 38 & $-58$ \\ \hline
coef. of $L_{\Q(\sqrt{-30})}(\phi, s)$ & $-2$ & 3 & 5 & 0 & $-2$ & $-14$ & 26 & 0 & 14 & $-38$ & $-58$ \\ \hline
\end{tabular}

\bigskip 
\end{center}

The final results yields
\begin{align}
L_{\Q(\sqrt{-30})}(\phi,3) &= 
L\left(f_{120}\otimes\left(\frac{-3}{\cdot}\right) ,3\right),\nonumber\\
\m(P_{18}) &= 
\frac{6\sqrt{120}}{\pi^3}L\left(f_{120}\otimes\left(\frac{-3}{\cdot}\right) ,3\right) +\frac{14}{5} d_3.
\label{eqn: P18 mahler L}
\end{align}

\section{Relating $L(\TT(Y),s)$ to a newform}

The main tool for this section is the following result from~\cite{Sch08}.

\begin{thm} [Sch\"utt] The following classification of singular $K3$-surfaces over $\Q$ are equivalent. 
\begin{itemize}
\item By the discriminant $d$ of the transcendental lattice of the surface up to square.

\item  By the discriminant $-d$ of the N\'eron-Severi lattice of the surface up to square.

\item By the associated newform up to twisting.

\item By the level of the associated newform up to square.

\item By the $CM$ field $\Q(\sqrt{-d})$ of the associated newform. 
\end{itemize}
\end{thm}

This theorem depends on Livn\'e's modularity theorem for singular $K3$-surfaces that predicts that $L(\TT(Y),s)$ is modular and that the corresponding modular form has weight 3.  

The first step in finding the corresponding modular form is to compute the first few coefficients $A_p$ from $L(\TT(Y),s)$; then the coefficients are compared to 
the tables that can be found in \cite{Sch08} in order to identify the corresponding CM newform. 
Tackling the first step requires the following result from~\cite{Be10}.

\begin{thm}[Bertin] 
Let $Y$ be an elliptic $K3$-surface defined over $\Q$ and rank  $r(Y)=0$. Then
\begin{equation}
A_p=-\sum_{s \in \mathbb{P}^1(\mathbb{F}_p)}a_p(s),
\label{eqn: Ap rank0}
\end{equation}
where
\[
a_p(s)=p+1-\# Y_s(\mathbb{F}_p).
\]

Now suppose that  $r(Y)=1$ and that there is an infinite section defined 
over $\Q(\sqrt{d})$. Then
\begin{equation}
A_p=-\sum_{s \in \mathbb{P}^1(\mathbb{F}_p)}a_p(s)-\left(\frac{d}{p}\right)p.
\label{eqn: Ap rank1}
\end{equation}
\end{thm}
Notice that the result stated in \cite{Be10} requires a {\em generator} of $\mathrm{MWL}(Y)$ to be defined over $\Q(\sqrt{d})$. But it is not hard to 
see that it suffices to find any element of infinite order to be defined over $\Q(\sqrt{d})$.

\subsection{Relating $L(\TT(Y_{6}),s)$ to a newform}
We know from Section~\ref{sec: lat and rank 6} that  $r(Y_6)=0$ and that $| \det \TT(Y_6) | = 24$, so we use equation~\eqref{eqn: Ap rank0}.  With the help of Pari/gp or Sage we compute several coefficients $A_p$ and compare them to the coefficients  of the newform of level 24 from 
Sch\"utt's table in~\cite{Sch08}.

\begin{center}

\begin{tabular}{|c || c| c|c| c| c| c| c| c| c| c| c|} \hline
$a_p$ &  5 & 7 & 11 & 13 & 17 & 19 & 23 & 29 & 31\\\hline \hline
newform of level 24  & $-2$ & $-10$ & 10 & 0 &0 & 0 & 0 & $-50$ & 38 \\ \hline
$A_p$ &  2 & $-10$ & $-10$ & 0 &0 & 0 & 0 & 50 & 38 \\ \hline
\end{tabular}

\end{center}

We see that
\begin{align*}
L(\TT(Y_6),3) &=L\left(f_{24}\otimes \left(\frac{-3}{\cdot}\right),3\right),\\
\intertext{and combining this with equation~\eqref{eqn: P6 mahler L} gives the final result}\\
\m(P_6) &= 
\frac{24\sqrt{6}}{\pi^3}L(\TT(Y_6),3).
\end{align*}

\subsection{Relating $L(\TT(Y_3),s)$ to a newform}
In this case, $r(Y_3)=1$ and the infinite section is defined over~$\Q$. We apply equation~\eqref{eqn: Ap rank1} to compute the $A_p$ values and  compare with the table from~\cite{Sch08} in order to obtain
\begin{align*}
L(\TT(Y_3),3) &= 
L\left(f_{15},3\right).\\
\intertext{Combining this with equation~\eqref{eqn: P3 mahler L} gives the final result}\\
\m(P_3) &= 
\frac{15\sqrt{15}}{2\pi^3}L(\TT(Y_3),3).
\end{align*}

\subsection{Relating $L(\TT(Y_{18}),s)$ to a newform} \label{Schutt-18}

In this case, $r(Y_{18})=1$ and the infinite section is defined over $\Q(\sqrt{-3})$.  We again apply equation~\eqref{eqn: Ap rank1} to compute the $A_p$ values and  compare with the table from \cite{Sch08} in order to obtain
\begin{align*}
L(\TT(Y_{18}),3) &= 
L\left(f_{120}\otimes\left(\frac{-3}{\cdot}\right) ,3\right).\\
\intertext{Combining this with equation~\eqref{eqn: P18 mahler L} gives the final result}\\
\m(P_{18}) &= 
\frac{120\sqrt{120}}{20\pi^3}L(\TT(Y_{18}),3) +\frac{14}{5} d_3.
\end{align*}

As a final note, we remark that one could have started the computations from this subsection without knowing that the infinite section is defined over $\Q(\sqrt{-3})$.  Computing several values of $A_p$ with equation~\eqref{eqn: Ap rank0} and comparing with the table from~\cite{Sch08} will reveal the necessary correction factor.
This allows one to \emph{predict} that the infinite section is defined over  $\Q(\sqrt{-3})$, and armed with this knowledge the infinite section is more easily computed
(see Section~\ref{sec:find inf sec}).

\section{Infinite section for $Y_{18}$}\label{sec:Inf Secs}
We now describe the computations used to  find an infinite section $p_\sigma$ for the elliptic surface given in equation~\eqref{eqn:Y18 model},   show that our $p_\sigma$ is a generator for the infinite part of the group of sections, and prove
that $h(p_\sigma)=10$. 

\subsection{Finding the infinite section}\label{sec:find inf sec}
 As noted above, we can predict that the infinite section 
is defined over $\Q(\sqrt{-3})$. Therefore, we  twist equation~\eqref{eqn:Y18 model} by $-3$ in order to get an elliptic surface with
the infinite section defined over $\Q$. 
We denote this twist $Y_{-3}$ (we drop the $Y_{18}$ notation in this case because there is no ambiguity).
Applying the general formula for a quadratic twist~\cite[Chapter~4]{Con}, we have
\[
Y_{-3}: 
y^2+(\sigma^2-18\sigma+1)xy
=x^3 + (-\sigma^4 + 36\sigma^3 - 329\sigma^2 + 90\sigma + 2)x^2 + 9\sigma(-\sigma + 18)x.
\]

For each $\sigma$, the fiber $Y_\sigma$ is a curve in $Y_{18}$ and the fiber $Y_{\sigma,-3}$ is a curve in $Y_{-3}$.  These curves satisfy the following exact sequence (see \cite{IR}, Proposition 20.5.4):
\[
0 \rightarrow 
Y_{\sigma,-3}(\Q) \rightarrow 
Y_\sigma\left(\Q(\sqrt{-3})\right) \stackrel{ \mathrm{Tr}_{\Q(\sqrt{-3})/\Q}}{\longrightarrow} 
Y_\sigma(\Q) \rightarrow 
Y_\sigma(\Q) / 2Y_\sigma(\Q) \rightarrow 0.
\]
More specifically, we have
\[
0 \rightarrow 
Y_{\sigma,-3}(\Q) \rightarrow 
Y_\sigma\left(\Q(\sqrt{-3})\right) \rightarrow \Z/6\Z \rightarrow \Z / 2\Z \rightarrow 0.
\]

A computation verifies that  a section for $Y_{-3}$ is given by  $p_{-3} = \left(x_{-3}(\sigma),y_{-3}(\sigma)\right)$ where
\begin{align*}
x_{-3}(\sigma)&=-\frac{2^4 3^6\sigma (\sigma-18) (\sigma-21)^2(\sigma+3)^2}{(\sigma-9)^2 (\sigma^2-21\sigma+72)^2 (\sigma^2-15\sigma+18)^2}, \text{ and} \\
y_{-3}(\sigma)&=-\frac{2^23^4\sigma (\sigma-18) (\sigma-21)(\sigma+3)}{(\sigma-9)^3 (\sigma^2-21\sigma+72)^3 (\sigma^2-15\sigma+18)^3} \\
& \qquad \cdot (\sigma^{10} - 108\sigma^9 + 4455\sigma^8 - 87822\sigma^7 + 771363\sigma^6 - 294840\sigma^5
- 44001711\sigma^4   \\
&\qquad\qquad \qquad  + 281168010\sigma^3 - 545848956\sigma^2 + 132322248\sigma + 128490624).
\end{align*}


The curve $Y_{\sigma,-3}$ has good reduction modulo 5 when $\sigma \equiv 1,2 \,(\mathrm{mod}\, 5)$. In those cases, one finds that  $Y_{\sigma,-3}(\F_5)$ has 6 elements and is
generated by the point $(3,1)$.  Hence the torsion of  $Y_{\sigma,-3}(\Q)$ injects into $\Z/6\Z$.
With the help of Pari/gp or Sage~\cite{PARI, St11}, it is easy to compute $[6]p_{-3}$ and see that the result is different from $O_{\sigma,-3}$. Therefore 
this point is not torsion.

Reversing the change of coordinates, one finds an infinite section $p_\sigma=(x(\sigma), y(\sigma))$ for the surface $Y_{18}$:
\begin{align}\label{infinitesec}
x(\sigma)&=\frac{2^4 3^5\sigma (\sigma-18) (\sigma-21)^2(\sigma+3)^2}{(\sigma-9)^2 (\sigma^2-21\sigma+72)^2 (\sigma^2-15\sigma+18)^2}, \text{ and}\\
y(\sigma)&=\frac{2^23^2 \sqrt{-3}\sigma(\sigma-21)(\sigma-18)(\sigma+3) }{(\sigma-9)^3 (\sigma^2-21\sigma+72)^3 (\sigma^2-15\sigma+18)^3}\left(\sigma^2+3(-6+\sqrt{-3})\sigma+9(5-3\sqrt{-3})\right)
\nonumber\\
& \qquad  
\cdot\left(\sigma^3+3(-9+\sqrt{-3})\sigma^2+9(19-6\sqrt{-3})\sigma+9(-9+11\sqrt{-3})\right) \nonumber \\
& \qquad \cdot \left(\sigma^5 +3(-15+4 \sqrt{-3})\sigma^4 +27(19-16\sqrt{-3})\sigma^3
 \right.\nonumber\\
 & \qquad \qquad\left. +81(9+52\sqrt{-3})\sigma^2  +162(-139-36\sqrt{-3})\sigma+5832(1-\sqrt{-3})\right).\nonumber
\end{align}




It is clear from these formulas that $p_\sigma$ and the zero section $[0:1:0]$ have simple intersections over $\sigma=9$, 
and over the distinct roots of  $(\sigma^2-21\sigma+72)$ and $ (\sigma^2-15\sigma+18)$.  Therefore $\overline{p_\sigma}\cdot\overline{O}=5$.
Applying equation~\eqref{heightdef}, we see that
\begin{eqnarray*}
h(p_\sigma) &=& 2 \chi(Y_{18})+ 2 \left(\overline{p_\sigma}\cdot \overline{O} \right)-\sum_{v}\mathrm{contr}_\nu(P)
= 2\cdot 2 + 2 \cdot 5 -\sum_{v}\mathrm{contr}_\nu(P) .\\
\end{eqnarray*}
From this, we have 
\begin{eqnarray*}
14\geq &h(p_\sigma)& \geq  14 -\frac{6 \cdot 6}{12}-\frac{1\cdot 1}{2}-\frac{1\cdot 1}{2}-\frac{1\cdot 2}{3}-\frac{1\cdot 2}{3}\\
14\geq  &h(p_\sigma)& \geq \frac{26}{3}.
\end{eqnarray*}
From the remarks following equation~\eqref{heightformula}, we know that the height of a generator must be either $5/2$ or $10$.
This means that $h(p_\sigma)=10$, since it must be a square multiple of the height of a generator. In Section~\ref{sec: height computation}, we show this fact directly by analyzing the intersection with the singular fibers. 


\subsection{Proof that $p_\sigma$ is a generator}\label{sec:generator}

Let $K = \Q(\sqrt{-3})(\sigma)$.  To prove that $p_\sigma$ is indeed a generator of the infinite section, we need to see that we cannot write $p_\sigma +k\rho_6= [2]P$ for any $P \in E(K)$ and $k=0, \dots, 5$. In fact, it suffices to prove that $p_\sigma +k\rho_6= [2]P$ has no solution $P \in E(K)$ for $k=0,3$. 
We will use the following theorem.
\begin{thm}[\cite{Con}, Proposition 1.7.5(b)]\label{thm:halving}
Let 
\[
E: y^2  = x(x^2 + ax +b)  
\]
 be an elliptic curve defined over a field $K$ with $\textup{char}\ K \neq 2$, and suppose $a^2 - 4b \not\in {K^*}^2$.  Let $Q = (x, y) \in E(K)$ with $x\neq 0$.  Then there exists $P \in E(K)$ such that $Q = [2]P $ iff \textup{(i)} $x \in {K^*}^2$, say $x = r^2$; and \textup{(ii)} one of $q_{\pm} = 2x + a \pm 2y/r \in {K^*}^2$.
\end{thm}

In order to  apply this result, we need to eliminate the term $xy$ from the Weierstrass equation~\eqref{eqn:Y18 model}, which we do by making the change $Y=y+\frac{(\sigma^2-18\sigma+1)x}{2}$. This gives
\[
Y^2= x\left(x^2 + \frac{\sigma^4 - 36\sigma^3 + 330\sigma^2 - 108\sigma - 3}{4}x + (-\sigma^2 + 18\sigma)\right).
 \]
From equation~\eqref{infinitesec}, we see that 
$x(\sigma)$ is not a square in $K$, hence there is no $P \in E(K)$ such that $p_\sigma = [2]P$.

Now write  $p_\sigma +3\rho_6=(x'(\sigma),Y'(\sigma))$.  A computation yields 
\[
x'(\sigma)
= -\frac{(\sigma-9)^2(\sigma^2-21\sigma+72)^2(\sigma^2-15\sigma+18)^2}{2^4\cdot 3^5 (\sigma-21)^2(\sigma+3)^2}
\]
which is a square in $K$, so take 
\[ r = 
\frac{(\sigma-9)(\sigma^2-21\sigma+72)(\sigma^2-15\sigma+18)}{2^2\cdot 3^2\sqrt{-3} (\sigma-21)(\sigma+3)}.
\] 
To compute $q_\pm$ as in Theorem~\ref{thm:halving}, we first find
\begin{align*}
Y'(\sigma)  &= 
\frac{ \sqrt{-3}(\sigma -9)  (\sigma^2 - 15\sigma + 18)(\sigma^2 - 21\sigma + 72) (\sigma^3 - 12\sigma^2 - 171\sigma + 1350)  }
{2^6 3^8(\sigma+3)^3 (\sigma -21)^3 }\\
 & \qquad \qquad\qquad\qquad \cdot 
(\sigma^3 - 42\sigma^2 +
369\sigma - 216)  (\sigma^4 - 36\sigma^3 + 351\sigma^2 - 486\sigma - 486).
\end{align*}
It is then a simple matter to compute
\begin{align*}
 q_+ & = -\frac{1}{2^2\cdot3^5}(\sigma-21)^2(\sigma+3)^2(\sigma^2 - 18\sigma + 9)\\
q_-&=-\frac{3^5(\sigma^2 - 18\sigma + 1)^3}{(\sigma-21)^2(\sigma+3)^2},
\end{align*}
and  neither of these are squares in $K$.




\subsection{Height computation}\label{sec: height computation}
In order to compute $h(p_\sigma)$, we need to study the intersection of $p_\sigma$ with the singular fibers,  since the correction terms in formula (\ref{heightdef}) are given by
\[
\mathrm{contr}_\nu(P)=\frac{j(m-j)}{m},
\]
when $P$ intersects the component $\Theta_{s,j}$ of the singular fiber over $s$ of type $I_m$.
We need the following theorem from~\cite{Ne64}:
\begin{thm}[N\'eron]
Let $E_s$ be an elliptic curve defined over $\mathbb{C}[s]$ 
given by a Weierstrass model, and denote by $v$ the $s$-adic valuation. Suppose that $E_0$ has a double point with distinct  tangents and 
$v(j(E_s))=-m<0$ \textup(this happens if and only if $E_0$ is singular of type $I_m$ in Kodaira's classification\textup). Then, for every integer $l>m/2$, there exists a Weierstrass model $\mathcal{E}_s$ deduced from $E_s$ by a transformation of the form
\begin{eqnarray*}
X&=&x+qz,\\
Y&=&y+ux+rz,\\
Z&=&z,\\
\end{eqnarray*}
with $q,r,u \in \mathbb{C}[s]$. A Weierstrass model $\mathcal{E}_s$ is given by 
\begin{equation}
Y^2Z+\lambda XYZ+\mu YZ^2=X^3 + \alpha X^2Z+\beta XZ^2+\gamma Z^3
\label{eqn: neron model eqn}
\end{equation}
with coeffcients satisfying 
\begin{equation}
v(\lambda^2+4 \alpha)=0, \quad 
v(\mu) \geq l,\quad  
v(\beta)\geq l,\quad  
v(\gamma)=m, \text{ and}\quad  
v(j(\mathcal{E}_s))=-m.
\label{eqn: neron model vals}
\end{equation}
\end{thm}

  We now follow the argument in~\cite{Be08English}, and refer the interested reader there for details.
A singular fiber of type $I_m$ over $s=0$ is composed of the nonsingular rational curves $\Theta_{0,0}, \Theta_{0,1},\dots, \Theta_{0,m-1}$.
If $m=2h$, the configuration of the these curves can be found in $\left(\PP^2\right)^h$, with a point $[X:Y:Z] \in Y_{18}$ over $s=0$ corresponding to the point 
\begin{equation}\label{eqn:neron over s=0}
[X:Y:Z^{(1)}]\times [X:Y:Z^{(2)}] \times \dots \times [X:Y:Z^{(h)}] \in (\mathbb{P}^2)^h,
\quad \text{ where }[X:Y:Z^{(i+1)}] = [X:Y:s Z^{(i)}].
\end{equation}
  So in particular, 
  \[
  [X:Y:Z^{(1)}] = [X:Y:sZ] \text{ and inductively } [X:Y:Z^{(h)}] = [X:Y:s^hZ].
  \]
If $[X:Y:Z]$ satisfies equation~\eqref{eqn: neron model eqn}, then $ [X:Y:Z^{(h)}]$ must satisfy the equation
\[
Y^2Z^{(h)}+\lambda XYZ^{(h)}+(\mu/s^h) Y(Z^{(h)})^2= s^h X^3 + \alpha X^2Z^{(h)}+ (\beta/s^h) X(Z^{(h)})^2+(\gamma/s^{2h}) (Z^{(h)})^3.
\]
Now, given the valuations in~\eqref{eqn: neron model vals} and the fact that $2h=m$, at $s=0$ this simplifies to
\begin{equation}
Y^2Z^{(h)}+\lambda_0 XYZ^{(h)}=  \alpha_0 X^2Z^{(h)}+ \gamma_m^0 (Z^{(h)})^3,
\label{eqn:conic}
\end{equation}
where the subscript~$0$ indicates evaluation at $s=0$, and  $\gamma_m^0 = (\gamma / s^{m})|_{s=0}$.

In fact, we can describe the components $\Theta_{0,i}$ exactly.  We give here only the fibers relevant in the sequel: 
\begin{align}
\Theta_{0,0} &= [X:Y:0]\times \dots \times [X:Y:0] \in (\mathbb{P}^2)^h, \text{ and}\nonumber\\
\Theta_{0,h} &=
[0:0:1]\times \dots \times [0:0:1] \times [X_0:Y_0:Z_0]\in (\mathbb{P}^2)^h,\label{eqn: component 6}
\end{align}
where $Z_0 \neq 0$ and $[X_0:Y_0:Z_0]$ is on the conic~\eqref{eqn:conic}.

\subsubsection{The fiber over $s=0$} This is a singularity of type $I_{12}$.  Let $\left(x'(\sigma), y'(\sigma)\right)$ represent the infinite section in equation~\eqref{infinitesec}.  The change of variables
\[
x(s) = s^4 x'(1/s), \qquad y(s) = s^6 y'(1/s), \qquad \sigma = 1/s
\]
  yields an infinite section for the Weierstrass model around $0$  given by the equation
\[y^2+(s^2-18s+1)xy=x^3 + s^2(-s^2 - 18s + 1)x^2 + (-s^6 + 18s^7)x.\]
A second change of variables 
\begin{align*}
x &=X+2s^6, &y&=Y-sX-2s^7-s^6
\end{align*} gives the $\mathcal{E}_s$ model
\begin{align*}
Y^2 + (s^2 &- 20s + 1)XY  + (2s^8 - 40s^7)Y =X^3 + (6s^6 - s^4 - 17s^3 - 18s^2 + s)X^2\\
&+ (12s^{12} - 4s^{10} - 68s^9 - 71s^8 + 2s^7)X  
+(8s^{18} - 4s^{16} - 68s^{15} - 70s^{14} - s^{12}).
\end{align*}
The same change of variables applied to the infinite section $\left( x(s), y(s) \right)$ yields
\[
\left( X(s) , Y(s) \right) = \left(s^6f_1(s),s^6g_1(s)\right)
\]
where $f_1(0)=-2$ and $g_1(0)=1$. So by equation~\eqref{eqn:neron over s=0} this corresponds to the point
\[
 [0:0:1]\times [0:0:1]\times [0:0:1]\times [0:0:1]\times [0:0:1]\times[-2:1:1]
\]
in the $\mathcal{E}_s$ model.  From~\eqref{eqn: component 6} we see that this point is on $\Theta_{0,6}$ because
$[-2:1:1]$ is on the conic
\[Y^2+XY+Z^2=0.\]

\subsubsection{The fiber over $s=\infty$} This is a singularity of type $I_2$, and the infinite section given in equation~\eqref{infinitesec} is for the model around infinity
given by the  Weierstrass equation  
\[y^2+(\sigma^2-18\sigma+1)xy=x(x-1)(x+\sigma^2-18\sigma)=x^3 + (\sigma^2 - 18\sigma - 1)x^2 + (-\sigma^2 + 18\sigma)x.\]
So we work with the singular fibers over $\sigma =0$ just as we did above with $s=0$.
The change of variables 
\begin{align*}
x &=\frac{X}{9}+12\sigma,
& 
y&=\frac{Y}{27}+\frac{X}{9}-6 \sigma
\end{align*}
 gives the $\mathcal{E}_\sigma$ model 
\begin{align*}
Y^2 + (3\sigma^2 - 54\sigma + 9)XY + &(324\sigma^3 - 5832\sigma^2)Y = X^3 + (324\sigma - 27)X^2\\
&+ (1458\sigma^3 + 8667\sigma^2)X + (157464\sigma^4 - 1583388\sigma^3 + 78732\sigma^2).
\end{align*}
The same change of variables applied to equation~\eqref{infinitesec} yields the infinite section
\[
\left( X(\sigma) , Y(\sigma) \right) = 
\left(
\sigma f_2(\sigma),\sigma g_2(\sigma)\right), \quad
\text{ where } f_2(0)=-\frac{1011}{8} \text{ and } g_2(0)=\frac{9099-1575\sqrt{-3}}{16}.
\]
 From~\eqref{eqn: component 6}, the corresponding point on  the $\mathcal{E}_\sigma$ model is
\[
  [0:0:1]\times  \left[-\frac{1011}{8}: \frac{9099-1575\sqrt{-3}}{16} :1\right],\]
which is on the component $\Theta_{\infty,1}$ since the second point is on the conic \[Y^2+9XY+27X^2-78732Z^2=0.\]

\subsubsection{The fiber over $s=\frac{1}{18}$} This is also a singularity of type $I_2$. We consider the change of variables
\begin{align}\label{eqn: coord change s=1/18}
X &=-y-(\sigma^2-18\sigma+1)x, &
Y&=y, &
Z&=x+(\sigma^2-18\sigma)z,
\end{align}
which takes the Weierstrass equation at infinity  to
\begin{equation}\label{eqn: alt infinite model}
(X+Y)(X+Z)(Y+Z)+(\sigma^2-18\sigma+1)XYZ=0.
\end{equation}
When $s = \frac{1}{18}$, we have $\sigma = 18$, and the equation is a product of two rational curves
\[
(X + Y + Z) (XY + XZ + YZ) = 0,
\]
so this is our N\'eron model.    The component $\Theta_{\frac{1}{18},0}$ is the one meeting the zero section, which is given by $[x:y:z] = [0:1:0]$.  From the change of coordinates in~\eqref{eqn: coord change s=1/18},  this corresponds to 
$
[X:Y:Z] = [-1:1:0]$.
  So we have
  \[
  \Theta_{\frac{1}{18},0}: X+Y+Z =0 \quad \text{ and } \quad \Theta_{\frac{1}{18},1}: XY + XZ + YZ=0.
  \]
Applying the change of coordinates in~\eqref{eqn: coord change s=1/18} to the infinite section in~\eqref{infinitesec}, one calculates 
\[ XY + XZ + YZ=-\frac{2^43^5(\sigma-18)\sigma(\sigma-21)^2(\sigma+3)^2}{(\sigma-9)^2 (\sigma^2-21\sigma+72)^2 (\sigma^2-15\sigma+18)^2},\]
which means that it cuts $\Theta_{\frac{1}{18},1}$. 



\subsubsection{The fibers over $s=\alpha_1$, $\beta_1$, $\alpha_2$, $\beta_2$} 
Recall that $\alpha_1$ and $\beta_1$ are the two distinct roots of  $s^2-18s+1=0$, and since $\sigma = 1/s$ they are also roots of  $\sigma^2-18\sigma+1$.  These fibers are of type $I_3$. We again use the change of coordinates in~\eqref{eqn: coord change s=1/18}.  From~\eqref{eqn: alt infinite model}, both fibers become a product of three rational curves
\[
(X+Y)(X+Z)(Y+Z) = 0.
\]
Again, the zero section is $[X:Y:Z] = [-1:1:0]$, which satisfies $X+Y=0$.  So we identify
  \[
  \Theta_{\alpha_1,0}: X+Y =0 \quad \text{ and } \quad \Theta_{\beta_1,0}: X+Y=0.
  \]
After the change of coordinates in~\eqref{eqn: coord change s=1/18}, the infinite section satisfies
\[X+Y=(\sigma^2-18\sigma+1)f_3(\sigma)\]
with $f_3(\sigma)$ a rational function not divisible by $(\sigma^2-18\sigma+1)$. 
Hence the infinite section cuts $\Theta_{\alpha_1,0}$ and~$\Theta_{\beta_1,0}$.

Finally, note that the fibers over $\alpha_2$ and $\beta_2$ are of type $I_1$, so we know that the infinite section cuts 
$\Theta_{\alpha_2,0}$ and~$\Theta_{\beta_2,0}$ because that is the only choice.

Recall from the discussion in section \ref{sec:Inf Secs}  that $\overline{p_\sigma}\cdot\overline{O}=5$. 
With these considerations, equation (\ref{heightdef}) tells us that
\[h(p_\sigma)=2\cdot 2 + 2 \cdot 5 -\frac{6\cdot 6}{12}-\frac{1\cdot 1}{2}-\frac{1\cdot 1}{2}=10,\]
which completes the proof.

\subsection*{Acknowledgements}  The authors would like to  thank the Banff International Research Station for sponsoring the second Women in Numbers workshop and for providing a productive and enjoyable environment for our initial work on this project.   We also thank Kiran Kedlaya, Joseph Silverman, and Bianca Viray for some helpful discussions.

\end{document}